\documentclass[preprint,3p,times]{elsarticle}

\usepackage{amssymb}
\usepackage{amsmath}
\usepackage[all,cmtip]{xy}
\usepackage{amsthm}
\usepackage{color}
\usepackage{enumitem}
\usepackage{mathtools}
\usepackage{tikz}
\usepackage{amsmath, mathrsfs}
\usepackage{appendix}
\usepackage[colorlinks=true]{hyperref}
\usepackage[hyphenbreaks]{breakurl}
\input diagxy

\theoremstyle{plain}
 \newtheorem{thm}{Theorem}[section]
 \newtheorem{lem}[thm]{Lemma}
 \newtheorem{prop}[thm]{Proposition}
 
\theoremstyle{definition}
 \newtheorem{defn}[thm]{Definition}
 \newtheorem{exmp}[thm]{Example}
 \newtheorem{rem}[thm]{Remark}
 
 \newtheorem{ques}[thm]{Question}

\DeclareMathAlphabet{\mathcal}{OMS}{cmsy}{m}{n}

\DeclareMathOperator{\ob}{ob}
\DeclareMathOperator{\Sub}{Sub}

\makeatletter
\def\ps@pprintTitle{%
\let\@oddhead\@empty
\let\@evenhead\@empty
\def\@oddfoot{\centerline{\thepage}}%
\let\@evenfoot\@oddfoot}
\makeatother

\def\oto{{\bfig\morphism<180,0>[\mkern-4mu`\mkern-4mu;]\place(86,0)[\circ]\efig}}

\newcommand{\da}{\downarrow}

\newcommand{\ra}{\rightarrow}

\newcommand{\lda}{\swarrow}
\newcommand{\rda}{\searrow}

\newcommand{\bv}{\bigvee}

\newcommand{\dv}{\dashv}

\newcommand{\nat}{\natural}

\newcommand{\dbv}{\displaystyle\bv}

\renewcommand{\phi}{\varphi}
\newcommand{\al}{\alpha}
\newcommand{\be}{\beta}

\newcommand{\ga}{\gamma}

\newcommand{\lam}{\lambda}
\newcommand{\si}{\sigma}

\newcommand{\Om}{\Omega}

\newcommand{\CE}{\mathcal{E}}

\newcommand{\CQ}{\mathcal{Q}}

\newcommand{\sD}{\mathsf{D}}

\newcommand{\sQ}{\mathsf{Q}}

\newcommand{\Fy}{\mathfrak{y}}

\newcommand{\Cat}{\mathbf{Cat}}

\newcommand{\Dist}{\mathbf{Dist}}

\newcommand{\Rel}{\mathbf{Rel}}
\newcommand{\Set}{\mathbf{Set}}
\newcommand{\Sh}{\mathbf{Sh}}

\newcommand{\QCat}{\CQ\text{-}\Cat}

\newcommand{\QCcCat}{\CQ\text{-}\mathbf{CcCat}}
\newcommand{\QCcSymCat}{\CQ\text{-}\mathbf{CcSymCat}}
\newcommand{\DQCcSymCat}{\DQ\text{-}\mathbf{CcSymCat}}

\newcommand{\QSet}{\sQ\text{-}\Set}
\newcommand{\OmSet}{\Om\text{-}\Set}

\newcommand{\QDist}{\CQ\text{-}\Dist}

\newcommand{\hf}{\widehat{f}}

\newcommand{\hX}{\widehat{X}}
\newcommand{\hY}{\widehat{Y}}
\newcommand{\hZ}{\widehat{Z}}
\newcommand{\hal}{\widehat{\alpha}}
\newcommand{\hbe}{\widehat{\beta}}
\newcommand{\hga}{\widehat{\gamma}}

\newcommand{\op}{{\rm op}}

\newcommand{\DQ}{\sD\sQ}

\newcommand{\DQCat}{\DQ\text{-}\Cat}

\newcommand{\with}{\mathrel{\&}}

\renewcommand{\leq}{\leqslant}

\newcommand{\sqleq}{\sqsubseteq}

\numberwithin{equation}{section}

\allowdisplaybreaks

\begin{document}

\begin{frontmatter}



\title{$\mathsf{Q}\text{-}\mathbf{Set}$ is not generally a topos}


\author{Xiao Hu}
\ead{huxiao97@stu.scu.edu.cn}

\author{Lili Shen\corref{cor}}
\ead{shenlili@scu.edu.cn}

\cortext[cor]{Corresponding author.}
\address{School of Mathematics, Sichuan University, Chengdu 610064, China}

\begin{abstract}
For a commutative, unital and divisible quantale $\mathsf{Q}$, it is shown that the category of $\mathsf{Q}$-sets is a topos if, and only if, $\mathsf{Q}$ is a frame.
\end{abstract}


\end{frontmatter}




\section{Introduction}

The categorical foundation plays a crucial role in the study of fuzzy sets \cite{Zadeh1965}. One of the most notable approaches is the theory of quantale-valued sets developed by H{\"o}hle and his collaborators \cite{Hoehle1991,Hoehle1992,Hoehle1995b,Hoehle1998,Hoehle2005,Hoehle2011a}, which extends the theory of frame-valued sets initiated by Higgs \cite{Higgs1970,Higgs1984} and Fourman--Scott \cite{Fourman1979}. 

Explicitly, from every frame $\Om$ we may construct a bicategory $\sD\Om$ \cite{Walters1981}, which is actually a \emph{quantaloid} \cite{Rosenthal1996,Stubbe2005,Stubbe2014}. Symmetric categories enriched in the quantaloid $\sD\Om$ are precisely $\Om$-sets, and morphisms between $\Om$-sets are exactly left adjoint $\sD\Om$-distributors. It is well known that the category
\[\OmSet\]
of $\Om$-sets and their morphisms is a topos (see \cite[Theorem 5.9 and Proposition 9.2]{Fourman1979}).

If we consider a unital and involutive \emph{quantale} \cite{Rosenthal1990,Mulvey1992} $\sQ$ as the table of truth values, the theory of $\sQ$-sets can be established in the same way. More specifically, we may construct a quantaloid $\DQ$ \cite{Hoehle2011a,Pu2012,Stubbe2014} (see Proposition \ref{DQ-def} in the case that $\sQ$ is commutative, unital and divisible), and define $\sQ$-sets as symmetric $\DQ$-categories \cite{Hoehle2011a,Pu2012} (see Definition \ref{Q-set-def} and Proposition \ref{Qset-as-sym-DQ-cat}). Thus, we obtain the category
\[\QSet\]
of $\sQ$-sets and left adjoint $\DQ$-distributors. It is now natural to ask:

\begin{ques}
Is $\QSet$ a topos?
\end{ques}

This question was first proposed by Pu--Zhang (see \cite[Question 7.2]{Pu2012}) in the case that $\sQ$ is a commutative, unital and divisible quantale (also known as \emph{GL-monoid} \cite{Hoehle1992}), where they conjectured that $\QSet$ is a topos only when $\sQ$ is a frame. The aim of this paper is to give an affirmative answer to their hypothesis (see Theorem \ref{QSet-topos-frame}):
\begin{itemize}
\item For a commutative, unital and divisible quantale $\sQ$, the category $\QSet$ is a topos if, and only if, $\sQ$ is a frame.
\end{itemize}
Therefore, the theory of topoi is unlikely to serve as the categorical foundation of $\sQ$-sets when the quantale $\sQ$ is non-idempotent.

\section{Divisible quantales}

A commutative and unital \emph{quantale} \cite{Mulvey1986,Rosenthal1990,Eklund2018} is a commutative monoid $(\sQ,\with,1)$ whose underlying set $\sQ$ is a complete lattice (with a top element $\top$ and a bottom element $\bot$), such that
\[p\with\Big(\bv_{i\in I}q_i\Big)=\bv_{i\in I}p\with q_i\]
for all $p,q_i\in\sQ$ $(i\in I)$. The right adjoint induced by the multiplication $\with$, denoted by $\ra$,
\[(p\with -)\dv(p\ra -)\colon\sQ\to\sQ,\]
satisfies
\[p\with q\leq r\iff p\leq q\ra r\]
for all $p,q,r\in\sQ$. We say that $\sQ$ is \emph{divisible} if it satisfies one of the following equivalent conditions:

\begin{prop} \label{divisible-def} (See \cite[Proposition 2.1]{Pu2012}.)
The following statements are equivalent:
\begin{enumerate}[label={\rm(\roman*)}]
\item \label{divisible-def:u-q} $u=q\with(q\ra u)$ whenever $u\leq q$ in $\sQ$.
\item \label{divisible-def:u-v-q} $v\with(q\ra u)=(q\ra v)\with u$ whenever $u,v\leq q$ in $\sQ$.
\item \label{divisible-def:u-p-q} $u=q\with p$ for some $p\in\sQ$ whenever $u\leq q$ in $\sQ$.
\item \label{divisible-def:p-q} $p\wedge q=p\with(p\ra q)$ for all $p,q\in\sQ$.
\end{enumerate}
In this case, the unit $1$ of the monoid $(\sQ,\with,1)$ must be the top element $\top$ of the complete lattice $\sQ$.
\end{prop}

\begin{rem}
As pointed out by an anonymous referee, a commutative and divisible quantale need not be unital; one of such examples is given in \cite[Exercise 2.2.1(15)]{Eklund2018}. Nevertheless, since we only consider unital quantales in this paper, we do not bother writing the non-unital version of Proposition \ref{divisible-def}.
\end{rem}

Throughout this paper, we always assume that
\[(\sQ,\with)\]
is a commutative, unital and divisible quantale (also \emph{GL-monoid} \cite{Hoehle1992}). 

\begin{exmp}
The following commutative, unital and divisible quantales are well known:
\begin{enumerate}[label=(\arabic*)]
\item If $\Om$ is a frame, then $(\Om,\wedge)$ is a commutative, unital and divisible quantale.
\item If $*$ is a \emph{continuous t-norm} \cite{Alsina2006,Klement2000,Klement2004b} on the unit interval $[0,1]$, then $([0,1],*)$ is a commutative, unital and divisible quantale \cite{Hoehle1992}.
\item The Lawvere quantale $([0,\infty],+)$ \cite{Lawvere1973} is commutative, unital and divisible.
\end{enumerate}
\end{exmp}

Let $q\in\sQ$. We say that $q$ is \emph{idempotent} if $q\with q=q$. In this case, it is easy to see that
\[p\with q=p\wedge q\]
for all $p\in\sQ$.

\begin{prop} \label{Q-underlying-frame} (See \cite[Theorem 5.2]{Hoehle1995a}.)
The underlying complete lattice of $\sQ$ is a frame. In particular, the idempotent elements of $\sQ$ form a subframe of the frame $(\sQ,\wedge)$.
\end{prop}

By Proposition \ref{divisible-def}\ref{divisible-def:u-v-q}, it makes sense to define
\begin{equation} \label{circ_q-def}
v\circ_q u\coloneqq v\with(q\ra u)=(q\ra v)\with u
\end{equation}
whenever $u,v\leq q$ in $\sQ$.

\begin{prop} \label{da-q-circ-q-divisible-quantale} (See \cite[Proposition 2.5]{Pu2012}.)
For every $q\in\sQ$, $(\da q,\circ_q)$ is a commutative, unital and divisible quantale, where $\da q=\{p\in\sQ\mid p\leq q\}$.
\end{prop}

Moreover, we say that $p\in\sQ$ is \emph{idempotent relative to} $q\in\sQ$ \cite{Pu2012} if $p$ is idempotent in the quantale $(\da q,\circ_q)$; that is, if
\begin{equation} \label{u-idem-to-q}
p\leq q\quad\text{and}\quad p\circ_q p=p\with(q\ra p)=p.
\end{equation}
Hence, by Propositions \ref{Q-underlying-frame} and \ref{da-q-circ-q-divisible-quantale}, the set
\begin{equation} \label{Cq-def}
C_q\coloneqq\{p\leq q\mid p\with(q\ra p)=p\}
\end{equation}
of idempotent elements relative to $q$ is a subframe of the underlying frame of the quantale $(\da q,\circ_q)$.

\begin{prop} \label{u-wedge-p}
Suppose that $q\in\sQ$ and $p\in C_q$. Then $p\wedge r=p\with(q\ra r)$ for all $r\in\sQ$.
\end{prop}

\begin{proof}
On one hand, since $1$ is the top element of $\sQ$ and $p\in C_q$, it is clear that
\[p\with(q\ra r)\leq p\with 1=p\quad\text{and}\quad p\with(q\ra r)\leq q\with(q\ra r)\leq r.\]
On the other hand,
\[p\wedge r=p\with(p\ra r)=p\with (q\ra p)\with(p\ra r)\leq p\with(q\ra r),\]
where the first equality follows from Proposition \ref{divisible-def}\ref{divisible-def:p-q}, and the second equality follows from $p\in C_q$ and \eqref{Cq-def}. 
\end{proof}

For any $p,q\in\sQ$, define
\begin{equation} \label{sqleq-def}
p\sqleq q\iff p\in C_q\iff p\leq q\ \text{and}\ p\with(q\ra p)=p.
\end{equation}

\begin{prop} \label{sqleq-po}
$\sqleq$ is a partial order on $\sQ$. Therefore, $C_q$ is the principal $\sqleq$-lower set generated by $q\in\sQ$.
\end{prop}

\begin{proof}
$\sqleq$ is clearly reflexive and antisymmetric. For the transitivity, suppose that $p\sqleq q$ and $q\sqleq r$. Then
\begin{align*}
p&=q\with(q\ra p) & (p\leq q\ \text{and Proposition \ref{divisible-def}\ref{divisible-def:u-q}})\\
&=q\with(r\ra q)\with(q\ra p) & (q\sqleq r\ \text{and \eqref{sqleq-def}})\\
&=(r\ra q)\with p & (p\leq q\ \text{and Proposition \ref{divisible-def}\ref{divisible-def:u-q}})\\
&=(r\ra q)\with(q\ra p)\with p & (p\sqleq q\ \text{and \eqref{sqleq-def}})\\
&\leq p\with(r\ra p).
\end{align*}
Since the reverse inequality is trivial, we conclude that $p\with(r\ra p)=p$. Hence $p\sqleq r$.
\end{proof}

\begin{lem} \label{pqr-sqleq}
If $p\leq q\leq r$ and $p\sqleq r$, then $p\sqleq q$.
\end{lem}

\begin{proof}
This is an immediate consequence of $p=p\with(r\ra p)\leq p\with(q\ra p)$.
\end{proof}

\begin{prop} \label{sqleq-join}
If $A\subseteq\sQ$ and $a\sqleq q$ for all $a\in A$, then $\bv A\sqleq q$. Therefore, if $A\subseteq\sQ$ has a $\sqleq$-upper bound, then the $\sqleq$-join $\sqcup A$ of $A\subseteq\sQ$ exists, and $\sqcup A=\bv A$. 
\end{prop}

\begin{proof}
Since $a\sqleq q$ for all $a\in A$, by \eqref{sqleq-def} we have $\bv A\leq q$ and
\[\mbox{$\bv A$}=\bv\limits_{a\in A}a\with(q\ra a)\leq\bv\limits_{a\in A}a\with(q\ra\mbox{$\bv A$})=(\mbox{$\bv A$})\with(q\ra\mbox{$\bv A$}).\]
Thus $\bv A\sqleq q$. Moreover, note that $a\sqleq\bv A$ for all $a\in A$ by applying Lemma \ref{pqr-sqleq} to $a\leq\bv A\leq q$ and $a\sqleq q$. Therefore, $\bv A$ is the $\sqleq$-join of $A$.
\end{proof}

\begin{prop} \label{sqleq-meet}
If $A\subseteq\sQ$ and $A\neq\varnothing$, then the $\sqleq$-meet $\sqcap A$ of $A$ exists.
\end{prop}

\begin{proof}
Since it is clear that $\bot\sqleq q$ for all $q\in\sQ$, $A$ has $\bot$ as its $\sqleq$-lower bound. We show that
\begin{equation} \label{sqcap-A-def}
\sqcap A\coloneqq\bv\{q\in\sQ\mid\forall a\in A\colon q\sqleq a\}
\end{equation}
is the $\sqleq$-meet of $A$. First, $\sqcap A$ is a $\sqleq$-lower bound of $A$ by Proposition \ref{sqleq-join}. Second, if $q\sqleq a$ for all $a\in A$, then $q\leq\sqcap A\leq a$ for all $a\in A$, and consequently $q\sqleq\sqcap A$ by Lemma \ref{pqr-sqleq}. 
\end{proof}






\section{Quantaloid-enriched categories and their Cauchy completeness}

A \emph{quantaloid} \cite{Rosenthal1996} $\CQ$ is a category whose hom-sets are complete lattices, such that the composition $\circ$ of $\CQ$-arrows preserves suprema on both sides, i.e.,
\[v\circ\Big(\bv_{i\in I} u_i\Big)=\bv_{i\in I}v\circ u_i\quad\text{and}\quad\Big(\bv_{i\in I} v_i\Big)\circ u=\bv_{i\in I}v_i\circ u\]
for all $\CQ$-arrows $u,u_i\colon p\to q$, $v,v_i\colon q\to r$ $(i\in I)$. The corresponding right adjoints induced by the composition maps 
\[(-\circ u)\dv(-\lda u)\colon \CQ(p,r)\to\CQ(q,r)\quad\text{and}\quad(v\circ -)\dv(v\rda -)\colon \CQ(p,r)\to\CQ(p,q)\] 
satisfy 
\[v\circ u\leq w\iff v\leq w\lda u\iff u\leq v\rda w\] 
for all $\CQ$-arrows $u\colon p\to q$, $v\colon q\to r$, $w\colon p\to r$. 

If a pair of $\CQ$-arrows $u\colon p\to q$ and $v\colon q\to p$ satisfy
\[1_p\leq v\circ u\quad\text{and}\quad u\circ v\leq 1_q,\]
we say that $u$ and $v$ form an adjunction in $\CQ$ and denote it by $u\dv v$, where $u$ is called a \emph{left adjoint} of $v$, and $v$ is a \emph{right adjoint} of $u$. In this case, it is easy to check that
\[v=u\rda 1_q\quad\text{and}\quad u=1_q\lda v.\]
Therefore, the right adjoint of a $\CQ$-arrow, when it exists, is necessarily unique. We define
\begin{equation} \label{u*-def}
u^*\coloneqq u\rda 1_q\colon q\to p
\end{equation}
for each $\CQ$-arrow $u\colon p\to q$.

\begin{lem} (See \cite[Proposition 2.3.4]{Heymans2010}.) \label{map-calculation}
If $u\colon p\to q$ is a left adjoint in $\CQ$, then 
\[v\circ u=v\lda u^*\quad\text{and}\quad u^*\circ w=u\rda w\]
for all $\CQ$-arrows $v\colon q\to r$, $w\colon r\to q$.
\end{lem}

We say that $\CQ$ is \emph{involutive} if $\CQ$ is equipped with an \emph{involution}
\[(-)^{\circ}\colon\CQ^{\op}\to\CQ,\]
which is a functor satisfying
\[q^{\circ}=q,\quad u^{\circ\circ}=u\quad\text{and}\quad\Big(\bigvee_{i\in I}u_i\Big)^{\circ}=\bigvee_{i\in I}u_i^{\circ}\]
for all $q\in\ob\CQ$ and $\CQ$-arrows $u,u_i\colon p\to q$ $(i\in I)$.

Given a small (i.e., $\ob\CQ$ is a set) and involutive quantaloid $\CQ$, a \emph{$\CQ$-category} (also \emph{category enriched in $\CQ$}) \cite{Stubbe2005} consists of a set $X$, a \emph{type} map $|\text{-}|:X\to\ob\CQ$ and a family of $\CQ$-arrows $\al(x,y)\in\CQ(|x|,|y|)$ $(x,y\in X)$, such that
\[1_{|x|}\leq\al(x,x)\quad\text{and}\quad\al(y,z)\circ\al(x,y)\leq\al(x,z)\]
for all $x,y,z\in X$. Note that $(X,\al)$ has an underlying (pre)order given by
\[x\leq y\iff |x|=|y|\ \text{and}\ 1_{|x|}\leq\al(x,y).\]
We say that $(X,\al)$ is \emph{separated} (or \emph{skeletal}) if the underlying order on $X$ is a partial order. Moreover, $(X,\al)$ is \emph{symmetric} \cite{Heymans2011} if
\begin{equation} \label{sym-Q-cat}
\al(x,y)=\al(y,x)^{\circ}
\end{equation}
for all $x,y\in X$.

A \emph{$\CQ$-functor} $f:(X,\al)\to(Y,\be)$ between $\CQ$-categories is a map $f:X\to Y$ such that
\begin{equation} \label{Q-functor-def}
|x|=|fx|\quad\text{and}\quad\al(x,y)\leq\be(fx,fy)
\end{equation}
for all $x,y\in X$. $\CQ$-categories and $\CQ$-functors constitute a category
\[\QCat.\]

A \emph{$\CQ$-distributor} $\phi\colon(X,\al)\oto(Y,\be)$ between $\CQ$-categories is a map that assigns to each pair $(x,y)\in X\times Y$ a $\CQ$-arrow $\phi(x,y)\in\CQ(|x|,|y|)$, such that\[\be(y,y')\circ\phi(x,y)\circ\al(x',x)\leq\phi(x',y')\]
for all $x,x'\in X$, $y,y'\in Y$. With the pointwise order inherited from $\CQ$, the category
\[\QDist\] 
of $\CQ$-categories and $\CQ$-distributors becomes a (large) quantaloid in which
\begin{align*}
&\psi\circ\phi\colon(X,\al)\oto(Z,\ga),\quad(\psi\circ\phi)(x,z)=\bigvee_{y\in Y}\psi(y,z)\circ\phi(x,y),\\
&\xi\lda\phi\colon(Y,\be)\oto(Z,\ga),\quad(\xi\lda\phi)(y,z)=\bigwedge_{x\in X}\xi(x,z)\lda\phi(x,y),\\
&\psi\rda\xi\colon(X,\al)\oto(Y,\be),\quad(\psi\rda\xi)(x,y)=\bigwedge_{z\in Z}\psi(y,z)\rda\xi(x,z)
\end{align*}
for all $\CQ$-distributors $\phi\colon(X,\al)\oto(Y,\be)$, $\psi\colon(Y,\be)\oto(Z,\ga)$, $\xi\colon(X,\al)\oto(Z,\ga)$; the identity $\CQ$-distributor on $(X,\al)$ is given by its hom $\al\colon(X,\al)\oto(X,\al)$. 

Adjoint $\CQ$-distributors are precisely adjunctions in the quantaloid $\QDist$. Each $\CQ$-functor $f\colon(X,\al)\to(Y,\be)$ induces a pair of adjoint $\CQ$-distributors $f_{\nat}\dv f_{\nat}^*$, given by
\[f_{\nat}\colon(X,\al)\oto(Y,\be),\quad f_{\nat}(x,y)=\be(fx,y)\quad\text{and}\quad f_{\nat}^*\colon(Y,\be)\oto(X,\al),\quad f_{\nat}^*(y,x)=\be(y,fx).\]
It is straightforward to verify the following lemma:

\begin{lem} \label{fully-faithful-graph}
A $\CQ$-functor $f\colon(X,\al)\to(Y,\be)$ is \emph{fully faithful} in the sense that
\[\al(x,x')=\be(fx,fx')\]
for all $x,x'\in X$ if, and only if, $f_{\nat}^*\circ f_{\nat}=\al$.
\end{lem}

For each $q\in\ob\CQ$, let $\{q\}$ denote the (necessarily symmetric) one-object $\CQ$-category whose only object has type $q$ and hom $1_q$. For every $\CQ$-category $(X,\al)$, it is easy to see that
\begin{equation} \label{al(x,-)-dv-al(-,x)}
\al(x,-)\dv\al(-,x)\colon(X,\al)\oto\{|x|\}
\end{equation}
for all $x\in X$. In fact, considering the $\CQ$-functor 
\begin{equation} \label{tilde-x-def}
x\colon\{|x|\}\to(X,\al)
\end{equation}
targeting at $x$, it is clear that $\al(x,-)=x_{\nat}$ and $\al(-,x)=x_{\nat}^*$.

A $\CQ$-category $(X,\al)$ is \emph{Cauchy complete} if every left adjoint $\CQ$-distributor
\[\mu\colon\{q\}\oto(X,\al)\]
is \emph{representable}; that is, there exists $x\in X$, called the \emph{representation} of $\mu$, such that $\mu=\al(x,-)$.

Let $(X,\al)$ be a $\CQ$-category. Define 
\begin{equation} \label{hX-def}
\hX\coloneqq\{\mu\colon\{q\}\oto(X,\al)\mid \mu\ \text{is a left adjoint in}\ \QDist,\ q\in\ob\CQ\}
\end{equation}
and
\begin{equation} \label{hal-def}
\hal(\mu,\lam)\coloneqq\lam\rda\mu=\lam^*\circ\mu=(\lam\rda\al)\circ\mu
\end{equation}
for all $\mu,\lam\in\hX$, where the last two equalities follow from Lemma \ref{map-calculation} and Equation \eqref{u*-def}, respectively. Then $(\hX,\hal)$ is a separated and Cauchy complete $\CQ$-category, called the \emph{Cauchy completion} of $(X,\al)$ (cf. \cite{Betti1982} and \cite[Proposition 4.12]{Pu2012}). There is a fully faithful $\CQ$-functor
\begin{equation} \label{Yoneda-embedding}
\Fy\colon(X,\al)\to(\hX,\hal),\quad \Fy x=\al(x,-),
\end{equation}
which is also injective when $(X,\al)$ is separated. Moreover, the assignment $(X,\al)\mapsto(\hX,\hal)$ is functorial; that is, there is a functor
\begin{equation} \label{hat-functor}
\widehat{(-)}\colon\QCat\to\QCcCat
\end{equation}
sending each $\CQ$-functor $f\colon(X,\al)\to(Y,\be)$ to
\begin{equation} \label{hf-mu-def}
\hf\colon(\hX,\hal)\to(\hY,\hbe),\quad \hf\mu=f_{\nat}\circ\mu,
\end{equation}
where $\QCcCat$ is the full reflective subcategory of $\QCat$ consisting of separated and Cauchy complete $\CQ$-categories (see \cite[Proposition 7.14]{Stubbe2005}). Hence, there is a bijection
\begin{equation} \label{hat-dv-inclusion}
\QCcCat((\hX,\hal),(Y,\be))\cong\QCat((X,\al),(Y,\be))
\end{equation}
natural in $(X,\al)\in\QCat$ and $(Y,\be)\in\QCcCat$. In particular, for a $\CQ$-functor $f\colon(X,\al)\to(Y,\be)$ with $(X,\al)\in\QCat$ and $(Y,\be)\in\QCcCat$,
\begin{itemize}
\item the transpose of $f$ is
\begin{equation} \label{overline-f-def}
\overline{f}\colon(\hX,\hal)\to(Y,\be),
\end{equation}
with $\overline{f}\mu$ being the representation of $\hf\mu$ (see \eqref{hf-mu-def}) for all $\mu\in\hX$, and
\item the triangle
\begin{equation} \label{fy=f}
\bfig
\qtriangle<700,400>[(X,\al)`(\hX,\hal)`(Y,\be);\Fy`f`\overline{f}]
\efig
\end{equation}
is commutative.
\end{itemize}

In this paper we are concerned with the full subcategory
\[\QCcSymCat\]
of $\QCat$, whose objects are separated, symmetric and Cauchy complete $\CQ$-categories. For later use, we point out that equalizers in $\QCcSymCat$ are formulated in the same way as in $\QCat$ (cf. \cite[Remark 2.4(3)]{Shen2015}):

\begin{prop} \label{QCcSymCat-equalizer}
The equalizer of $f,g\colon(X,\al)\to(Y,\be)$ in $\QCcSymCat$ is given by
\[E=\{x\in X\mid fx=gx\}\]
equipped with the restriction of the $\CQ$-category structure $\al$ on $E$. Hence, the embedding of the equalizer into $(X,\al)$ is a fully faithful $\CQ$-functor.
\end{prop}

\begin{proof}
Let $\si=\al|(E\times E)$. By \cite[Remark 2.4(3)]{Shen2015}, it suffices to show that $(E,\si)$ is Cauchy complete. Let $\mu\colon\{q\}\oto(E,\si)$ be a left adjoint $\CQ$-distributor. Considering the inclusion $\CQ$-functor $j\colon(E,\si)\ \to/^(->/(X,\al)$, we have a left adjoint $\CQ$-distributor
\[j_{\nat}\circ\mu\colon\{q\}\oto(X,\al).\]
Since $(X,\al)$ is Cauchy complete, there exists $x\in X$ such that $j_{\nat}\circ\mu=\al(x,-)$. It follows that
\[\be(fx,-)=f_{\nat}\circ\al(x,-)=f_{\nat}\circ j_{\nat}\circ\mu=g_{\nat}\circ j_{\nat}\circ\mu=g_{\nat}\circ\al(x,-)=\be(gx,-),\]
where the third equality follows from $fj=gj$. Thus, the separatedness of $(Y,\be)$ forces $fx=gx$, and consequently $x\in E$. Therefore,
\[\mu=\si\circ\mu=j_{\nat}^*\circ j_{\nat}\circ\mu=j_{\nat}^*\circ\al(x,-)=\al(x,j-)=\si(x,-),\]
where the second equality follows from Lemma \ref{fully-faithful-graph}, as desired.
\end{proof}

\section{Quantale-valued sets as enriched categories}

We are now ready to introduce the key notion of this paper:

\begin{defn} (See \cite{Hoehle1992,Hoehle1995b,Hoehle2011a,Pu2012}.) \label{Q-set-def}
A \emph{$\sQ$-set} a set $X$ equipped with a map
\[\al\colon X\times X\to\sQ,\]
such that
\begin{enumerate}[label=(S\arabic*)]
\item \label{Q-set:d} $\al(x,y)\leq\al(x,x)\wedge\al(y,y)$,
\item \label{Q-set:s} $\al(x,y)=\al(y,x)$,
\item \label{Q-set:t} $\al(y,z)\with(\al(y,y)\ra\al(x,y))\leq\al(x,z)$
\end{enumerate}
for all $x,y,z\in X$. 
\end{defn}

In order to exhibit $\sQ$-sets as enriched categories, we need the following quantaloid constructed from $\sQ$:

\begin{prop} \label{DQ-def} (See \cite{Hoehle2011a,Pu2012,Stubbe2014}.)
The following data define an involutive quantaloid $\DQ$:
\begin{itemize}
\item $\ob\DQ=\sQ$;
\item $\DQ(p,q)=\{u\in\sQ\mid u\leq p\wedge q\}$;
\item the composite of $u\in\DQ(p,q)$ and $v\in\DQ(q,r)$ is given by $v\circ_q u$ (see \eqref{circ_q-def});
\item the identity $\DQ$-arrow of $\DQ(q,q)$ is $q$ itself;
\item each hom-set $\DQ(p,q)$ is equipped with the order inherited from $\sQ$;
\item the involution on $\DQ$ is the identity functor $1_{\DQ}\colon(\DQ)^{\op}=\DQ\to\DQ$.
\end{itemize}
\end{prop}

From the definition we see that a $\DQ$-category consists of a set $X$, a map $|\text{-}|\colon X\to\sQ$ and a map $\al\colon X\times X\to\sQ$ such that
\[\al(x,y)\leq|x|\wedge|y|,\quad |x|\leq\al(x,x)\quad\text{and}\quad\al(y,z)\with(|y|\ra\al(x,y))\leq\al(x,z)\]
for all $x,y,z\in X$. Note that the first and the second inequalities above force
\begin{equation} \label{alxx=|x|}
\al(x,x)=|x|
\end{equation}
for all $x\in X$. Thus, a $\DQ$-category is exactly given by a map $\al\colon X\times X\to\sQ$ satisfying \ref{Q-set:d} and \ref{Q-set:t} of Definition \ref{Q-set-def}, and a $\sQ$-set is precisely a $\DQ$-category satisfying \ref{Q-set:s}. Therefore:

\begin{prop} (See \cite[Proposition 6.3]{Hoehle2011a}.) \label{Qset-as-sym-DQ-cat}
A $\sQ$-set is precisely a symmetric $\DQ$-category.
\end{prop}

We denote by
\[\QSet\]
the category whose objects are $\sQ$-sets, and whose morphisms are left adjoint $\DQ$-distributors.

\begin{rem} \label{OmSet}
If $\sQ=\Om$ is a frame, then we obtain the well-known category $\OmSet$ as considered in \cite{Fourman1979} and \cite[Sections 2.8 and 2.9]{Borceux1994c}. Since $\OmSet$ is equivalent to the category $\Sh(\Om)$ of sheaves on $\Om$ (cf. \cite[Theorem 5.9]{Fourman1979} and \cite[Theorem 2.9.8]{Borceux1994c}), it follows that $\OmSet$ is a topos (see \cite[Proposition 9.2]{Fourman1979} and \cite[Example 5.2.3]{Borceux1994c}).
\end{rem}

\begin{rem}
Since the Cauchy completion preserves the symmetry of a $\DQ$-category (see Proposition \ref{Sym-Cc-Sym} infra), it can be deduced from \eqref{hal-def} that the right adjoint $\phi^*$ of a left adjoint $\DQ$-distributor $\phi\colon(X,\al)\oto(Y,\be)$ between $\sQ$-sets has the form 
\[\phi^*(y,x)=\phi(x,y)\] 
for all $x\in X$, $y\in Y$ (cf. \cite[Proposition 5.5]{Pu2012}). Therefore, it is straightforward to check that a $\DQ$-distributor $\phi\colon(X,\al)\oto(Y,\be)$ is a left adjoint if, and only if,
\begin{enumerate}[label=(\arabic*)]
\item $\al(x,x)=\dbv\limits_{y\in Y}\phi(x,y)\with(\be(y,y)\ra\phi(x,y))$,
\item $\phi(x,y)\with(\al(x,x)\ra\phi(x,y'))\leq\be(y,y')$
\end{enumerate}
for all $x\in X$, $y,y'\in Y$.
%
%
%
\end{rem}

Given a $\sQ$-set $(X,\al)$, a left adjoint $\DQ$-distributor $\mu\colon\{q\}\oto(X,\al)$ is usually called a \emph{singleton} \cite{Hoehle2011a,Pu2012}. In elementary words:

\begin{defn} \label{singleton-def} (See \cite[Definition 6.4]{Hoehle2011a} and \cite[Definition 5.3]{Pu2012}.) 
A \emph{singleton} on a $\sQ$-set $(X,\al)$ is precisely a map $\mu\colon X\to\sQ$ such that
\begin{enumerate}[label={\rm(ss\arabic*)}]
\item \label{singleton-def:d} $\mu(x)\leq\al(x,x)$,
\item \label{singleton-def:dist} $\mu(x)\with(\al(x,x)\ra\al(x,y))\leq\mu(y)$,
\item \label{singleton-def:a1} $|\mu|=\dbv\limits_{x\in X}\mu(x)\with(\al(x,x)\ra\mu(x))$,
\item \label{singleton-def:a2} $\mu(x)\with(|\mu|\ra\mu(y))\leq\al(x,y)$
\end{enumerate}
for all $x,y\in X$, where $|\mu|=\dbv\limits_{x\in X}\mu(x)$. 
\end{defn}

\begin{rem} \label{singleton-representation}
A singleton $\mu$ on $(X,\al)$ described by Definition \ref{singleton-def} is actually a left adjoint $\DQ$-distributor
\[\mu\colon\{|\mu|\}\oto(X,\al).\]
If $\mu$ is representable, then there exists $x\in X$ such that 
\[\mu(y)=\al(x,y)\in\DQ(|\mu|,|y|)\]
for all $y\in X$, which necessarily forces the type of $x$ (cf. \eqref{alxx=|x|}) to be
\begin{equation} \label{|x|=|mu|}
\al(x,x)=|\mu|=\bv\limits_{x\in X}\mu(x).
\end{equation}
\end{rem}

Note that $\QSet$ and $\DQCcSymCat$ are exactly the categories $\mathbf{S}M\text{-}\mathbf{Mod}_{\mathrm{ad}}$ and $\mathbf{C}M\text{-}\Set$ in \cite{Pu2012}, respectively. Hence:

\begin{prop} \label{QSet-DQCcSymCat} (See \cite[Proposition 5.7(2)]{Pu2012}.)
The category $\QSet$ is equivalent to the full subcategory
\[\DQCcSymCat\]
of $\DQCat$ whose objects are separated and Cauchy complete $\sQ$-sets.
\end{prop}


From the construction of $\DQ$ (Proposition \ref{DQ-def}) and \cite[Remark 2.4(1)]{Shen2015} we see that the terminal object of the category $\DQCat$ is given by 
\[(\sQ,\wedge),\]
whose underlying set is $\sQ$ equipped with the identity map as its type map, and whose hom-arrows are given by $p\wedge q$ for all $p,q\in\sQ$. The following proposition tells us that $(\sQ,\wedge)$ is also the terminal object of $\DQCcSymCat$:

\begin{prop} \label{Q-wedge-Cc} (See \cite[Example 4.8]{Pu2012}.)
$(\sQ,\wedge)$ is a separated and Cauchy complete $\sQ$-set.
\end{prop}

The following proposition guarantees that the Cauchy completion of a $\sQ$-set is also a $\sQ$-set:

\begin{prop} \label{Sym-Cc-Sym} (See \cite[Theorem 6.6]{Hoehle2011a} and \cite[Proposition 5.4]{Pu2012}.)
The Cauchy completion $(\hX,\hal)$ of a symmetric $\DQ$-category $(X,\al)$ is also symmetric.
\end{prop}

By Definition \ref{singleton-def} it is straightforward to check that
\begin{equation} \label{p-sqleq-q-singleton}
p\sqleq q\iff p\ \text{is a singleton on the one-element}\ \sQ\text{-set}\ \{q\}.
\end{equation}
Thus, the following characterization of the set $C_q$ (see \eqref{Cq-def}) can be deduced from \eqref{hal-def} in combination with Propositions \ref{u-wedge-p} and \ref{DQ-def}:

\begin{prop} \label{Cq-Cc} (See \cite[Example 4.13]{Pu2012}.)
For each $q\in\sQ$, the Cauchy completion of the one-element $\sQ$-set $\{q\}$ is precisely $C_q$ equipped with the $\sQ$-set structure inherited from $(\sQ,\wedge)$.
\end{prop}

For $p,q\in\sQ$, since a $\DQ$-functor must preserve types, a $\DQ$-functor from $(C_p,\wedge)$ to $(C_q,\wedge)$, whenever it exists, is necessarily neutral on every element of $C_p$. Hence, there is at most one $\DQ$-functor from $(C_p,\wedge)$ to $(C_q,\wedge)$, which can only be the inclusion map and exists only when $C_p\subseteq C_q$. Correspondingly, the order $\sqleq$ may be described as follows:

\begin{prop} \label{sqleq-Cc}
Let $p,q\in\sQ$. Then $p\sqleq q$ if, and only if, there exists a (unique) $\DQ$-functor from $(C_p,\wedge)$ to $(C_q,\wedge)$.
\end{prop}

\begin{proof}
If $p\sqleq q$, then by Proposition \ref{sqleq-po}, every $r\in C_p$ belongs to $C_q$; that is, $C_p\subseteq C_q$. Thus, the inclusion map is the unique $\DQ$-functor from $(C_p,\wedge)$ to $(C_q,\wedge)$. Conversely, if the inclusion map is a $\DQ$-functor from $(C_p,\wedge)$ to $(C_q,\wedge)$, then $p\in C_q$, which means that $p\sqleq q$.
\end{proof}

The following proposition reveals that in a separated and Cauchy complete $\sQ$-set $(X,\al)$, an element $x\in X$ has a unique ``restriction'' on each $p\sqleq\al(x,x)$:

\begin{prop} \label{x_p-def}
Let $(X,\al)$ be a separated and Cauchy complete $\sQ$-set. For each $x\in X$ with $\al(x,x)=q$, if $p\sqleq q$, then there exists a unique element $x_p\in X$ such that
\begin{equation} \label{x_p-property}
\al(x_p,x_p)=p\quad\text{and}\quad\al(x_p,y)=\al(x,y)\wedge p=\al(x,y)\circ_q p
\end{equation}
for all $y\in X$.
\end{prop}

\begin{proof}
For each $y\in Y$, since $\al(x,y)\leq\al(x,x)=q$ and $p$ is idempotent in the quantale $(\da q,\circ_q)$ (see Proposition \ref{da-q-circ-q-divisible-quantale} and \eqref{u-idem-to-q}), we have
\[\al(x,y)\wedge p=\al(x,y)\circ_q p.\]
Thus, $\al(x,-)\wedge p$ is the composite of singletons
\[\bfig
\morphism<500,0>[\{p\}`\{q\};p]
\morphism(500,0)<600,0>[\{q\}`(X,\al),;\al(x,-)]
\place(250,0)[\circ]
\place(800,0)[\circ]
\efig\]
which is again a singleton on $(X,\al)$. Hence, the existence and uniqueness of $x_p\in X$ satisfying \eqref{x_p-property} follow from the Cauchy completeness and separatedness of $(X,\al)$.
\end{proof}

\begin{prop} \label{Cq-to-Y}
Let $q\in\sQ$, and let $(Y,\be)$ be a separated and Cauchy complete $\sQ$-set. Then every $\DQ$-functor from $(C_q,\wedge)$ to $(Y,\be)$ is uniquely determined by an element $y\in Y$ with $\be(y,y)=q$, which is exactly
\begin{equation} \label{overline-y-def}
\overline{y}\colon(C_q,\wedge)\to(Y,\be),\quad p\mapsto y_p
\end{equation}
where $y_p$ is the unique element determined by \eqref{x_p-property}.
\end{prop}

\begin{proof}
Since a $\DQ$-functor from $\{q\}$ to $(Y,\be)$ is essentially an element $y\in Y$ with $\be(y,y)=q$, the conclusion is an immediate consequence of \eqref{hat-dv-inclusion} and \eqref{overline-f-def} together with Propositions \ref{Cq-Cc} and \ref{x_p-def}.
\end{proof}

\begin{prop} \label{Cq-to-hatY}
Let $q\in\sQ$, and let $(Y,\be)$ be a $\sQ$-set. Then every $y\in Y$ with $\be(y,y)=q$ induces a $\DQ$-functor
\begin{equation} \label{hat-y-def}
\widehat{y}\colon(C_q,\wedge)\to(\hY,\hbe),\quad p\mapsto\be(y,-)\wedge p.
\end{equation}
\end{prop}

\begin{proof}
Note that $\widehat{y}$ is obtained by applying \eqref{hf-mu-def} to the $\DQ$-functor
\[y\colon\{q\}\to(Y,\be)\]
targeting at $y$. Indeed, by \eqref{hf-mu-def} we have
\[(\widehat{y}p)(y')=y_{\nat}(y')\circ_q p=\be(y,y')\circ_q p=\be(y,y')\wedge p\]
for all $p\in C_q$, $y'\in Y$, where the last equality holds because $\be(y,y')\leq\be(y,y)=q$ and $p$ is idempotent in the quantale $(\da q,\circ_q)$ (see Proposition \ref{da-q-circ-q-divisible-quantale}).
\end{proof}

\section{$\QSet$ is a topos only if $\sQ$ is a frame}

Let us recall the following facts in a topos:

\begin{prop} (see \cite[Proposition IV.1.2]{MacLane1992}.) \label{topos-monic}
Every monomorphism in a topos is an equalizer.
\end{prop}
 
\begin{prop} (see \cite[Proposition IV.6.3]{MacLane1992} and \cite[Proposition 5.10.2]{Borceux1994c}.) \label{topos-intersection-union}
Let $\CE$ be a topos. For each $X\in\ob\CE$, the partially ordered set $\Sub X$ is a lattice, in which the intersection $U\cap V$ and the union $w\colon U\cup V\to/>->/X$ of subobjects $u\colon U\to/>->/X$ and $v\colon V\to/>->/X$ are formulated as follows:
\begin{enumerate}[label={\rm(\arabic*)}]
\item \label{topos-intersection-union:intersection} $U\cap V$ given by the pullback below;
\[\bfig
\square/>->`>->`>->`>->/<600,400>[U\cap V`V`U`X;``v`u]
\efig\]
\item \label{topos-intersection-union:union} $U\cup V$ is given by the pushout on the left below, and $w\colon U\cup V\to/>->/X$ is determined by the universal property of the pushout, as on the right below. In particular, $w$ is a monomorphism.
\[\bfig
\square/>->`>->`>->`>->/<600,400>[U\cap V`V`U`U\cup V;```]
\square(1500,0)/>->`>->`>->`>->/<600,400>[U\cap V`V`U`U\cup V;```]
\morphism(1500,0)|b|/{@{>->}@/^-1em/}/<1100,-300>[U`X;u]
\morphism(2100,400)|a|/{@{>->}@/^1em/}/<500,-700>[V`X;v]
\morphism(2100,0)|b|/>-->/<500,-300>[U\cup V`X;w]
\efig\]
\end{enumerate}
\end{prop}

\begin{lem} \label{C_q-subobject}
Suppose that every monomorphism in $\DQCcSymCat$ is an equalizer. Then for each $q\in\sQ$, every subobject of $(C_q,\wedge)$ in $\DQCcSymCat$ is of the form $(C_p,\wedge)$ for some $p\sqleq q$.
\end{lem}

\begin{proof}
By Propositions \ref{QCcSymCat-equalizer} and \ref{topos-monic}, we may assume that a subobject of $(C_q,\wedge)$ in $\DQCcSymCat$ is a subset $S\subseteq C_q$ equipped with the inherited $\DQ$-category structure $\wedge$. We show that $S=C_p$ for some $p\sqleq q$.

First, $S$ has a $\sqleq$-maximum element $p$. Note that $S$ has a $\sqleq$-upper bound $q$, and thus, by Proposition \ref{sqleq-join}, $p:=\bv S$ is the $\sqleq$-join of $S$. We claim that $p\in S$. Indeed, by checking the conditions of Definition \ref{singleton-def} we see that
\[\mu\colon S\to\sQ,\quad \mu(s)=s\]
is a singleton on $(S,\wedge)$. Thus, the Cauchy completeness of $(S,\wedge)$ guarantees the existence of a representation $p'\in S$ of $\mu$. Then, it follows from Remark \ref{singleton-representation} (see \eqref{|x|=|mu|}) that 
\[p=\bv S=\bv\limits_{s\in S}\mu(s)=|\mu|=p'\wedge p'=p'\in S.\]



Second, $S=C_p$. It suffices to show that $r\in S$ for all $r\sqleq p$, which follows by applying Proposition \ref{x_p-def} to $(S,\wedge)$. Indeed, 
for each $r\sqleq p$ there exists a unique element $p_r\in S$ such that the first equality of \eqref{x_p-property} holds; that is,
\[r=p_r\wedge p_r=p_r\in S,\]
as desired.
%
\end{proof}

\begin{lem} \label{p-sqcap-q=p-wedge-q}
Suppose that $\DQCcSymCat$ is a topos. Then
\[p\sqcap q=p\wedge q\]
for all $p,q\in\sQ$.
\end{lem}

\begin{proof}
By Propositions \ref{QCcSymCat-equalizer}, \ref{Cq-Cc} and \ref{topos-monic}, $(C_p,\wedge)$ and $(C_q,\wedge)$ are subobjects of $(\sQ,\wedge)$ in $\DQCcSymCat$. Our strategy is to explore their union in $\Sub(\sQ,\wedge)$ under the guidance of Proposition \ref{topos-intersection-union}.


{\bf Step 1.} The intersection of $(C_p,\wedge)$ and $(C_q,\wedge)$ in $\Sub(\sQ,\wedge)$ is $(C_{p\sqcap q},\wedge)$. By Lemma \ref{C_q-subobject}, this intersection must be of the form $(C_r,\wedge)$ for some $r\sqleq p,q$. We claim that $r=p\sqcap q$. To this end, it suffices to show that the inner square of the diagram
\begin{equation} \label{p-sqcap-q-pullback}
\bfig
\pullback|brrb|<700,400>[(C_{p\sqcap q},\wedge)`(C_q,\wedge)`(C_p,\wedge)`(\sQ,\wedge);```]
|abb|/{@{->}@/^1em/}`-->`{@{->}@/^-1em/}/<600,400>[(X,\al);g`h`f]
\efig
\end{equation}
is a pullback. In fact, the commutativity of the inner square of \eqref{p-sqcap-q-pullback} is an immediate consequence of Propositions \ref{sqleq-meet} and \ref{sqleq-Cc}, in which every arrow is the inclusion $\DQ$-functor. Moreover, if the outer quadrilateral of \eqref{p-sqcap-q-pullback} is commutative, then 
\[hx\coloneqq fx=gx\in C_p\cap C_q\] 
for all $x\in X$. Thus, it follows from \eqref{sqleq-def} and Proposition \ref{sqleq-meet} that $hx\sqleq p\sqcap q$, i.e., $hx\in C_{p\sqcap q}$. Therefore, $h$ is the unique map making the triangles in \eqref{p-sqcap-q-pullback} commutative, and its $\DQ$-functoriality follows immediately from that of $f$ and $g$.

{\bf Step 2.} The union of $(C_p,\wedge)$ and $(C_q,\wedge)$ in $\Sub(\sQ,\wedge)$ is given by the Cauchy completion of the $\sQ$-set $(Z,\ga)$ with
\[Z=\{p,q\}\quad\text{and}\quad\ga(p,p)=p,\quad\ga(q,q)=q,\quad\ga(p,q)=\ga(q,p)=p\sqcap q.\]
Since $\DQCcSymCat$ is a topos, by Proposition \ref{topos-intersection-union}\ref{topos-intersection-union:union} it is sufficient to verify that the inner square of the diagram
\begin{equation} \label{p-union-q-pushout}
\bfig
\square<700,400>[(C_{p\sqcap q},\wedge)`(C_q,\wedge)`(C_p,\wedge)`(\hZ,\hga);``\widehat{q}`\widehat{p}]
\morphism|b|/{@{->}@/^-1em/}/<1300,-400>[(C_p,\wedge)`(Y,\be);\overline{y}]
\morphism(700,400)|a|/{@{->}@/^1em/}/<600,-800>[(C_q,\wedge)`(Y,\be);\overline{z}]
\morphism(700,0)|b|/-->/<600,-400>[(\hZ,\hga)`(Y,\be);\overline{h}]
\efig
\end{equation}
is a pushout, where $\widehat{p}$ and $\widehat{q}$ are the $\DQ$-functors defined by applying \eqref{hf-mu-def} to the inclusion $\DQ$-functors $\{p\}\ \to/^(->/(Z,\ga)$ and $\{q\}\ \to/^(->/(Z,\ga)$ targeting at $p$ and $q$, respectively.

First, the inner square of \eqref{p-union-q-pushout} is commutative. Suppose that $r\in C_{p\sqcap q}$. Then $r\sqleq p\sqcap q\sqleq p$. Thus, it follows from Proposition \ref{Cq-to-hatY} that
\[(\widehat{p}r)(p)=\ga(p,p)\wedge r=p\wedge r=r=(p\sqcap q)\wedge r=\ga(p,q)\wedge r=(\widehat{p}r)(q).\]
Similarly, $(\widehat{q}r)(p)=(\widehat{q}r)(q)=r$. Hence $\widehat{p}r=\widehat{q}r$, as desired.

Second, suppose that the outer quadrilateral of \eqref{p-union-q-pushout} is commutative, where $\overline{y}$ and $\overline{z}$ are the $\DQ$-functors determined by $y\in Y$ with $\be(y,y)=p$ and $z\in Y$ with $\be(z,z)=q$, respectively (see Proposition \ref{Cq-to-Y}). Then, Proposition \ref{Cq-to-Y} implies that
\begin{equation}\label{yr-equal-zr}
y_r=\overline{y}r=\overline{z}r=z_r
\end{equation}
for all $r\in C_{p\sqcap q}$; in particular,
\begin{equation} \label{y-z-p-sqcap-q}
y_{p\sqcap q}=z_{p\sqcap q}.
\end{equation}
Note that
\begin{equation} \label{h-Z-ga-Y-be}
h\colon(Z,\ga)\to(Y,\be),\quad hp=y\quad\text{and}\quad hq=z
\end{equation}
is a $\DQ$-functor, because
\[\ga(p,p)=p=\be(y,y)=\be(hp,hp),\quad\ga(q,q)=q=\be(z,z)=\be(hq,hq)\]
and
\begin{align*}
\ga(p,q)&=p\sqcap q\\
&=q\wedge(p\sqcap q)\\
&=\be(z,z)\wedge(p\sqcap q)&(\be(z,z)=q)\\
&=\be(z_{p\sqcap q},z)&(\text{Proposition \ref{x_p-def}})\\
&=\be(y_{p\sqcap q},z)&(\text{Equation \eqref{y-z-p-sqcap-q}})\\
&=\be(y,z)\wedge(p\sqcap q)&(\text{Proposition \ref{x_p-def}})\\
&\leq\be(y,z)\\
&=\be(hp,hq).
\end{align*}
We assert that the transpose 
\[\overline{h}\colon(\hZ,\hga)\to(Y,\be)\]  
of $h$ (see \eqref{hat-dv-inclusion} and \eqref{overline-f-def}) is the unique $\DQ$-functor such that the two triangles in \eqref{p-union-q-pushout} are commutative.
Indeed, by applying the functor $\widehat{(-)}$ to the commutative triangles
\[\bfig
\Vtrianglepair|aalrr|/^(->`<-^)`->`->`->/<700,400>[\{p\}`(Z,\ga)`\{q\}`(Y,\be);p`q`y`h`z]
\efig\]
we obtain that 
\[\widehat{h}\widehat{p}=\widehat{y}\quad\text{and}\quad\widehat{h}\widehat{q}=\widehat{z}.\] Then, by the definition of $\overline{(-)}$ (see \eqref{overline-f-def}) we conclude that 
\[\overline{h}\widehat{p}=\overline{y}\quad\text{and}\quad\overline{h}\widehat{q}=\overline{z}.\]For the uniqueness of $\overline{h}$, suppose that $\overline{k}\colon(\hZ,\hga)\to(Y,\be)$ is the transpose of $k\colon(Z,\ga)\to(Y,\be)$ (see \eqref{hat-dv-inclusion} and \eqref{overline-f-def}) and satisfies $\overline{k}\widehat{p}=\overline{y}$ and $\overline{k}\widehat{q}=\overline{z}$. Then it follows soon from the adjoint property (i.e., the natural bijection \eqref{hat-dv-inclusion}) that $kp=y$ and $kp=z$. Hence $k=h$.

{\bf Step 3.} Since $(\hZ,\hga)$ is the union of $(C_p,\wedge)$ and $(C_q,\wedge)$ in $\Sub(\sQ,\wedge)$, it is also a subobject of $(\sQ,\wedge)$, with the monomorphism
\[\overline{h}\colon(\hZ,\hga)\to(\sQ,\wedge)\]
given by setting $(Y,\be)=(\sQ,\wedge)$ in \eqref{p-union-q-pushout} (see Proposition \ref{topos-intersection-union}). Note that $\overline{h}$ is the transpose of the $\DQ$-functor (see \eqref{h-Z-ga-Y-be})
\[h\colon(Z,\ga)\to(Q,\wedge),\quad hp=p\quad\text{and}\quad hq=q,\]
and Propositions \ref{QCcSymCat-equalizer} and \ref{topos-monic} guarantee that $\overline{h}$ is fully faithful. Thus, by the fully faithfulness of \eqref{Yoneda-embedding} and the commutativity of the triangle \eqref{fy=f}, the composite $\DQ$-functor
\[h=\Big((Z,\ga)\to^{\Fy}(\hZ,\hga)\to^{\overline{h}}(\sQ,\wedge)\Big)\]
is fully faithful. Hence,
\[p\sqcap q=\ga(p,q)=hp\wedge hq=p\wedge q,\]
which completes the proof.
\end{proof}

\begin{prop} \label{main}
$\DQCcSymCat$ is a topos if, and only if, $\sQ$ is a frame.
\end{prop}

\begin{proof}
The ``if'' part is a direct consequence of Remark \ref{OmSet} and Proposition \ref{QSet-DQCcSymCat}. For the ``only if'' part, suppose that $\DQCcSymCat$ is a topos. Setting $p=1$ in Lemma \ref{p-sqcap-q=p-wedge-q} we obtain that $q\sqcap 1=q\wedge 1=q$, and consequently $q\sqleq 1$ for all $q\in\sQ$. By \eqref{sqleq-def}, this means that
\[q\with q=q\with(1\ra q)=q\]
for all $q\in\sQ$; that is, every element of $\sQ$ is idempotent. Hence $\sQ$ is a frame.
\end{proof}

Therefore, the main result of this paper is an immediate consequence of Propositions \ref{QSet-DQCcSymCat} and \ref{main}:

\begin{thm} \label{QSet-topos-frame}
$\QSet$ is a topos if, and only if, $\sQ$ is a frame.
\end{thm}

\begin{rem}
Our proof of Theorem \ref{QSet-topos-frame} relies on the fact that every topos satisfies the properties stated in Propositions \ref{topos-monic} and \ref{topos-intersection-union}. However, as pointed out by an anonymous referee, being an \emph{adhesive category} \cite{Lack2005} is sufficient to guarantee the validity of Propositions \ref{topos-monic} and \ref{topos-intersection-union} (see \cite[Lemma 4.8 and Theorem 5.1]{Lack2005}). Note that every topos is adhesive (see \cite[Proposition 3.7]{Lack2005}), and there are adhesive categories which are not topoi (see \cite[Example 3.8]{Lack2005}). Thus, Theorem \ref{QSet-topos-frame} can be slightly strengthened to:
\begin{itemize}
\item $\QSet$ is an adhesive category if, and only if, $\sQ$ is a frame.
\end{itemize}
\end{rem}

\begin{rem}
An anonymous referee suggests an alternative approach to proving Theorem \ref{QSet-topos-frame} by contradiction. Suppose that $\DQCcSymCat$ is a topos, but $\sQ$ is not a frame. Then there exists $p\in\sQ$ such that $p\not\sqleq 1$. 
\begin{itemize}
\item First, it can be calculated that $(\widehat{W},\wedge)$ with
\[\widehat{W}=\{r\in\sQ\mid r=((p\wedge r)\with(p\ra r))\vee(r\with r)\}\]
is the union of $(C_p,\wedge)$ and $(C_1,\wedge)$ in $\Sub(\sQ,\wedge)$, which is actually the Cauchy completion of $(W,\wedge)$ with $W=\{p,1\}$.
\item Second, it can be checked that the square
\[\bfig
\square<700,400>[(C_p\cap C_1,\wedge)`(C_1,\wedge)`(C_p,\wedge)`(\widehat{W},\wedge);```]
\efig\]
is not a pushout in $\DQCcSymCat$. 
\item Hence, it follows from Proposition \ref{topos-intersection-union} that $\DQCcSymCat$ cannot be a topos, contradicting the hypothesis.
\end{itemize}
%
%
\end{rem}

\begin{rem}
There is another possible direction towards the proof of Theorem \ref{QSet-topos-frame}. Explicitly, let
\[\DQ\text{-}\mathbf{SymDist}\]
denote the quantaloid of symmetric $\DQ$-categories and $\DQ$-distributors. Then it can be deduced from \cite[Theorem 4.7]{Heymans2012} that $\DQ\text{-}\mathbf{SymDist}$ is equivalent to the quantaloid $\Rel(\CE)$ of internal relations in a topos $\CE$ if, and only if, $\sQ$ is a frame. Hence, if we can figure out the relationship between $\DQ\text{-}\mathbf{SymDist}$ and the category $\Rel(\QSet)$ of internal relations in $\QSet$, then Theorem \ref{QSet-topos-frame} could possibly be established in light of \cite[Theorem 4.7]{Heymans2012}.
\end{rem}

\section*{Acknowledgement}

The authors acknowledge the support of National Natural Science Foundation of China (No. 12071319) and the Fundamental Research Funds for the Central Universities (No. 2021SCUNL202). The authors would like to thank Hongliang Lai, Isar Stubbe, Dexue Zhang, Jian Zhang and Jie Zhang for helpful discussions. The authors are also grateful for helpful remarks received from the anonymous referees which improve the presentation of this paper.






\end{document}